\documentclass[11pt]{amsart}

\usepackage{geometry}
\usepackage[french]{babel}

\usepackage{times}
\usepackage[scaled=0.92]{helvet}

\newtheorem*{thm}{Th\'eor\`eme}

\usepackage{amsfonts}
\usepackage{amssymb}
\renewcommand{\tilde}{\widetilde}

\DeclareMathOperator{\Aut}{\mathrm{Aut}}
\DeclareMathOperator{\Hom}{\mathrm{Hom}}
\DeclareMathOperator{\maks}{\mathfrak{m}}
\DeclareMathOperator{\kate}{\mathsf{Art}_{\mathit k}}
\DeclareMathOperator{\psl}{\![\![\!}
\DeclareMathOperator{\psr}{\!]\!]}

\DeclareMathOperator{\Set}{\mathsf{Ens}}

\def\og{\leavevmode\raise.3ex\hbox{$\scriptscriptstyle\langle\!\langle$~}}
\def\fg{\leavevmode\raise.3ex\hbox{~$\!\scriptscriptstyle\,\rangle\!\rangle$}}

\begin{document}

\date{\today\ (version 1.0)} 
\title[Un anneau de d\'eformation universel en conducteur sup\'erieur]{Un anneau de d\'eformation universel en conducteur sup\'erieur 
}

\author[J.~Byszewski, G.~Cornelissen, F.~Kato]{Jakub Byszewski, Gunther Cornelissen et Fumiharu Kato}
\address{Mathematisch Instituut, Universiteit Utrecht, Postbus 80.010, 3508 TA Utrecht, Nederland}
\email{j.j.byszewski@uu.nl}
\address{Mathematisch Instituut, Universiteit Utrecht, Postbus 80.010, 3508 TA Utrecht, Nederland}
\email{g.cornelissen@uu.nl}
\address{Department of Mathematics,  Faculty of Sciences, University of Kyoto, Kyoto 606-8502, Japan}
\email{kato@math.kyoto-u.ac.jp}
\subjclass[2000]{}

\begin{abstract} 
\noindent Soit $k$ un corps parfait de caract\'eristique $5$. Nous d\'emontrons que l'anneau versel de l'action d'un \'el\'ement d'ordre $5$ et de conducteur de Hasse $2$ comme automorphisme d'un anneau de s\'eries formelles $k\psl t \psr$ calcul\'e par Bertin et M\'ezard, est en fait \emph{uni}versel. C'est le premier exemple d'anneau non-trivial de d\'eformation universel  en conducteur sup\'erieur. 

\medskip 

\begin{center} {\large A universal deformation ring in higher conductor} \end{center}

\medskip 

\noindent {\sc Abstract.} Let $k$ denote a perfect field of characteristic $5$. We show that the versal deformation ring of an element of order $5$ and Hasse conductor $2$ as automorphism of a ring of formal power series $k\psl t \psr$, computed by Bertin and M\'ezard, is in fact \emph{uni}versal. This provides the first example of a non trivial universal deformation ring in higher conductor. 

 \end{abstract}

\maketitle

\section{Introduction} 

On ne connait que quelques anneaux de d\'eformation formelle pour les actions de groupes finis par automorphismes (continus) d'un anneau de s\'eries formelles en caract\'eristique positive: pour les actions de conducteur de Hasse $1$ (dites \og faiblement ramifi\'ees \fg , \cite{BM}, \cite{CK}, \cite{CM}, \cite{BC}), et pour l'action d'un \'element d'ordre $5$ et conducteur $2$ en caract\'eristique $5$ (\cite{BM}). N\'eanmoins, ces actions sont interessantes; par exemple,  le principe \og local-global \fg\,  de Henrio, Green-Matignon et Bertin-M\'ezard (\cite{Henrio}, \cite{GM}, \cite{BM}) implique que les questions de d\'eformations et de rel\`evements des actions de groupes sur les courbes alg\'ebriques se r\'eduisent \`a des questions similaires pour les actions par automorphismes de $k\psl t \psr$. Dans cette note, nous d\'emontrons que l'anneau versel de d\'eformation d'une action d'ordre $5$ et conducteur $2$ en caract\'eristique $5$ est universel. Ceci donne un premier exemple non trivial d'un anneau universel pour une action non faiblement ramifi\'ee. Les foncteurs de d\'eformations locales admettant une infinitude d'automorphismes infinit\'esimales, ce r\'esultat peut \^etre une surprise. La d\'emonstration d'universalit\'e pour les actions faiblement ramifi\'ees dans \cite{BC} utilise des m\'ethodes qui ne se g\'en\'eralisent pas aux cas de conducteur sup\'erieur et donc  ne donne pas de  raisons \'evidentes pour \'etablir l'universalit\'e dans une situation plus g\'en\'erale. Pour ces raisons, nous  esp\'erons que notre r\'esultat, qui n'est au fond qu'un calcul, serve \`a illuminer la question de l'universalit\'e en conducteur sup\'erieur.

Pourquoi s'interesser \`a l'universalit\'e? Remarquons que, pour les questions g\'eom\'etriques, la versalit\'e des anneaux de d\'eformation suffit souvent pour obtenir les r\'esultats d\'esir\'es. Par contre, en th\'eorie des nombres, on a besoin d'\'etablir l'universalit\'e. Si les d\'eformations de repr\'esentations lin\'eaires d'un groupe sont comprises (\og Lemme de Schur \fg , cf.\ \cite{Mazur}), les repr\'e\-sen\-ta\-tions par automorphismes des s\'eries formelles restent souvent plus myst\'erieuses. Pourtant  l'universalit\'e est parfois essentielle.  Par exemple, dans  \cite[Remark 3.3]{BC})  pour calculer l'anneau \emph{versel} de la d\'eformation d'une action faiblement ramifi\'ee de $\mathbf{Z}/p \oplus \mathbf{Z}/p$, on utilise l'\emph{universalit\'e} de l'action de $\mathbf{Z}/p$. Un autre exemple: pour \'etablir des r\'esultats de \og d\'evissage \fg \mbox{ } pour les foncteurs de d\'eformations locales il est parfois n\'ecessaire de supposer que les anneaux de d\'eformations sont universels \cite[Theorem 6.4.7]{Bthese}, \cite{B}. 

\section{Universalit\'e} Soit $k$ un corps parfait de caract\'eristique 5 et soit $W(k)$ l'anneau des vecteurs de Witt de $k$. Tout automorphisme de l'anneau de s\'eries formelles $k\psl t \psr$ d'ordre 5 et de conducteur (de Hasse) $\mathrm{ord}_t(\frac{\sigma(t)}{t}-1)$ \'egal \`a 2 est conjug\'e \`a : \begin{equation} \label{enp} \sigma \, : \, t \mapsto \frac{t}{\sqrt{t^2+1}}. \end{equation} (Voir \cite[4.2.1]{BM}).) Bertin et M\'ezard (\cite{BM}) ont \'etudi\'e le foncteur de d\'eformations formelles infinit\'esimales de $\sigma$. Soit $\kate$ la cat\'egorie des $W(k)$-alg\`ebres locales artiniens de corps r\'esiduel $k$. Pour un anneau $A$ de $\kate$, on appelle \emph{rel\`evement} de $\sigma$ \`a $A$ une s\'erie formelle $\tilde{\sigma}(t) \in \Aut_A A\psl t \psr$ d'ordre 5 (pour la composition) et telle  que $\tilde{\sigma}(t) \equiv {\sigma(t)} \mbox{ mod } \maks_A$. On dit que deux rel\`evements $\tilde{\sigma}_1$ et $\tilde{\sigma}_2$ de $\sigma$ \`a $A$ sont \'equivalents s'ils sont conjugu\'es par un \'el\'ement $\xi$ de $\Aut_A A\psl t \psr$ tel que $\xi(t) \equiv t \mbox{ mod } \maks_A$. On d\'efinit le foncteur de d\'eformations formelles infinit\'esimales de $\sigma$
$$D_\sigma : \kate \rightarrow \Set,$$
qui associe \`a un anneau $A$ de $\kate$ l'ensemble des classes d'\'equivalence de rel\`evements de $\sigma$ \`a $A$. 

Rappelons quelques notions de th\'eorie  infinit\'esimale des foncteurs (voir aussi \cite{schlessinger}). On dit qu'un foncteur $D \colon \kate \to \Set$ est \emph{pro-repr\'esentable}, s'il existe un anneau local noetherien complet $R$ de corps residuel $k$ tel que $D$ est isomorphe \`a un foncteur $\Hom(R,\cdot)$, qui associe \`a un anneau $A$ de $\kate$ l'ensemble des homomorphismes locaux de $W(k)$-alg\`ebres de $R$ dans $A$. Si un tel anneau existe, on appelle $R$ l'\emph{anneau universel} de $\sigma$. Parfois, il n'est facile que d'\'etablir un condition plus faible : on dit qu'un foncteur $D \colon \kate \to \Set$ admet un \emph{anneau versel} $R$, si $D(k)$ est r\'eduit \`a un seul \'el\'ement et s'il existe un foncteur pro-repr\'esentable $F=\Hom(R,\cdot)$ et un morphisme lisse $\varphi\colon F \to D$, c.-\`a-d., tel que pour tout morphisme surjectif $A' \to A$ dans $\kate$, l'application induite $$F(A')\to D(A') \times_{D(A)} F(A)$$ est surjective ; de plus, on demande que l'application $F(k[\varepsilon]/\varepsilon^2) \to D(k[\varepsilon]/\varepsilon^2)$ soit bijective. Lorsqu'il en est ainsi, on appelle $\varphi$ un \emph{morphisme versel}. L'anneau $R$ et le morphisme $\varphi$ sont alors unique (\`a un isomorphisme non-unique pr\`es), et le foncteur $D$ est pro-repr\'esentable si et seulement si le morphisme versel est un isomorphisme. On appelle l'\emph{action verselle} une classe de conjugaison d'une s\'erie formelle $\tilde{\sigma} \in \Aut_R R\psl t \psr$ telle que l'image de $\tilde{\sigma}$ dans $D(R/\maks_R^n)$ est \'egale \`a l'image par le morphisme versel de projection canonique $R \to R/\mathfrak{m}^n$. 

Bertin et M\'ezard ont d\'emontr\'e (\cite[Th\'eor\`eme 4.2.8]{BM}) que $D_{\sigma}$ admet un anneau versel 
$$ R_\sigma = W(k)[y]/\langle 1+y+y^2+y^3+y^4 \rangle $$ avec l'action verselle 
$$ \sigma_y \, : \, t \mapsto \frac{t}{\sqrt{t^2+y}}. $$ Le cas de caract\'eristique 5 et conducteur 2 est le seul cas du conducteur $>1 $ dans lequel l'anneau versel d'une action d'un groupe cyclique a \'et\'e compl\`etement d\'etermin\'e.

\begin{thm}
L'anneau $R_\sigma$ est \emph{uni}versel. 
\end{thm}

\begin{proof}
Soit $A$ un anneau de la cat\'egorie $\kate$, $\maks_A$ son id\'eal maximal. Nous prouvons que le morphisme versel est un isomorphisme. Supposons que $y_1,y_2 \in A$ satisfont $y_1^4+y_1^3+y_1^2+y_1+1=y_2^4+y_2^3+y_2^2+y_2+1=0$, $y_1 \equiv y_2 \equiv 1 \mbox{ mod }\maks_A$ et qu'il existe une s\'erie formelle $g(t)\in A\psl t \psr$ inversible (pour la composition) tel que $g(t) \equiv t \mbox{ mod } \maks_A$ et
$$ \sigma_{y_2} \circ g = g \circ \sigma_{y_1}. $$
Pour d\'emontrer l'universalit\'e, il suffit de montrer que $y_1=y_2$ --- la seule complication \'etant que l'on consid\`ere un anneau $A \in \kate$ arbitraire (le r\'esultat est trivial sur un corps). Si l'on pose $g(t)=a_0+a_1t+a_2t^2+O(t^3)$, on observe que $a_0 \in \maks_A$ et $a_1 \in A^*$. Les termes d'ordre 0 et 1 du d\'eveloppement en $t$ de l'\'egalit\'e 
\begin{equation}\label{surykatka} \frac{g(t)}{\sqrt{g(t)^2+y_1}} = g\left(\frac{t}{\sqrt{t^2+y_2}}\right). \end{equation}
donnent les formules :
\begin{equation} \label{zeroth}
a_0 = \frac{a_0}{\sqrt{a_0^2+y_1}} 
\end{equation}
et
\begin{equation} \label{wiewiorka} a_1 \frac{1}{\sqrt{a_0^2+y_1}}-\frac{a_0}{(a_0^2+y_1)^{3/2}}a_0a_1 = a_1 \frac{1}{\sqrt{y_2}}. \end{equation}
Dans (\ref{wiewiorka}), on a $a_1$ inversible et on peut utiliser (\ref{zeroth}) pour simplifier (\ref{wiewiorka}), ainsi \begin{equation} \label{first} \frac{1}{\sqrt{a_0^2+y_1}}-\frac{1}{\sqrt{y_2}}=a_0^2. \end{equation}
Les termes d'ordre trois dans (\ref{surykatka}) donnent
\begin{eqnarray*}
 & & \frac{{a_0}  {a_1}^2}{\left({a_0   }^2+{y_1}\right)^{3/2}}  
 -\frac{{a_0}}{\sqrt   {{a_0}^2+{y_1}}}
   \left(\frac{{a_0}^2   {a_1}^2}{\left({a_0   }^2+{y_1}\right)^2}+\frac{1}{2}   \left(\frac{{a_0}^2
   {a_1}^2}{\left({a_0   }^2+{y_1}\right)^2}-\frac{{a_1}^2+2 {a_0}   {a_2}}{{a_0}^2+{y_1}}\right)\right) \\ & & = \frac{{a_2}}{\sqrt{{  a_0}^2+{y_1}}}-\frac{{a_2}}{{y_2}}.
\end{eqnarray*}
En utilisant (\ref{zeroth}) et (\ref{first}),  
\begin{equation} \label{second}
a_2 \frac{1}{\sqrt{y_2}} \left( \frac{1}{\sqrt{y_2}} -1 \right) = \frac32 (a_0^2-1) a_0 a_1^2.
\end{equation}
Comme $\frac 32 (a_0^2-1)a_1^2$ est inversible, on a $$a_0 \in \left( \frac{1}{\sqrt{y_2}} -1 \right) A.$$ Multipliant (\ref{first}) par $a_0$ et utilisant (\ref{zeroth}), on obtient $$a_0 \left( \frac{1}{\sqrt{y_2}} -1 \right) = a_0^3,$$ et par suite $a_0^2 \in a_0^3 A$. Comme l'anneau $A$ est local, cela entra\^ine $a_0^2=0$ ; la conclusion r\'esulte de (\ref{first}). 
\end{proof}

\bibliographystyle{amsplain}\providecommand{\bysame}{\leavevmode\hbox to3em{\hrulefill}\thinspace}
\providecommand{\MR}{\relax\ifhmode\unskip\space\fi MR }
\providecommand{\MRhref}[2]{%
  \href{http://www.ams.org/mathscinet-getitem?mr=#1}{#2}
}
\providecommand{\href}[2]{#2}

\end{document}